\documentclass{amsart}
\title{An analysis of a war-like card game}
\author{Boris Alexeev \and Jacob Tsimerman}
\address{Department of Mathematics\\
Princeton University\\
Fine Hall, Washington Road\\
Princeton, NJ 08544-1000}
\email{balexeev@math.princeton.edu, jtsimerm@math.princeton.edu}
\date{December 2, 2009}

\makeatletter
\@namedef{subjclassname@2010}{%
  \textup{2010} Mathematics Subject Classification}
\makeatother


\usepackage{fullpage, setspace}
\usepackage{pstricks}

\renewcommand{\phi}      {\ensuremath{\varphi    }}
\newcommand{\halfopen}[2]{\ensuremath{\left(#1,#2\right]}}

\newcommand{\strong}[1]{\textbf{#1}}

\newtheorem*{maintheorem}         {Main Theorem}
\newtheorem*{mainlemma}           {Main Lemma}
\newtheorem*{monotonicitylemma}   {Monotonicity Lemma}
\theoremstyle{definition}
\newtheorem*{definition*}         {Definition}
\newtheorem*{problem*}            {Problem}

\begin{document}

\begin{abstract}
  In his book ``Mathematical Mind-Benders'', Peter Winkler poses the
  following open problem, originally due to the first author: ``[In the
    game Peer Pressure,] two players are dealt some number of cards,
  initially \emph{face up}, each card carrying a different integer.  In
  each round, the players simultaneously play a card; the higher card is
  discarded and the lower card passed to the other player.  The player
  who runs out of cards loses.  As the number of cards dealt becomes
  larger, what is the limiting probability that one of the players will
  have a winning strategy?''

  We show that the answer to this question is zero, as Winkler
  suspected.  Moreover, assume the cards are dealt so that one player
  receives $r \ge 1$ cards for every one card of the other.  Then if $r
  < \phi = \frac{1+\sqrt{5}}{2}$, the limiting probability that either
  player has a winning strategy is still zero, while if $r > \phi$, it
  is one.
\end{abstract}
\maketitle

\section*{Introduction}

The card game ``Peer Pressure'', a variant of ``War'', is played with a
deck of $n$ cards, each carrying a distinct integer.  The cards are
initially dealt randomly to two players, with either exactly $n/2$ cards
per player or each card randomly going to one of the players.  In each
round (``battle'') of the game, the players simultaneously play a card.
The player holding the higher card wins the round and receives the lower
card; however, the higher card is permanently discarded from the game.
The player who runs out of cards loses.  We assume that both players
know the original deck and thus are aware of the contents of both
players hands at all times.

Recall that we say a player has a \emph{winning strategy} if she may
announce her strategy beforehand and still always win against any
strategy from her opponent.  A winning strategy may in general be a
mixed (probabilistic) strategy, but if one exists, there also exists a
pure (deterministic) winning strategy.

In Peer Pressure, if there are four or fewer cards, then one of the
players has a winning strategy.  In particular, if one player has more
cards, then she wins; if the players have an equal number of cards, the
player with the highest card wins.  However, if there are five cards,
there is one position where neither player has a winning strategy:
$1,2,4$ versus $3,5$.

Suppose our two players are named Alice and Bob and have $a$ and $b$
cards respectively.  We prove the following lemma that helps classify
when a player has a winning strategy:

\begin{mainlemma}
  Let $\phi = \frac{1+\sqrt{5}}{2} \approx 1.61803$ be the golden ratio,
  which notably satisfies $1 + \phi = \phi^2$.
  \begin{itemize}
  \item
    If Alice has more than $\phi$ times as many cards as Bob (that is,
    $a > \phi b$), then Alice has a winning strategy.
  \item
    If Alice has more than $1/\phi$ times as many cards as Bob (that is,
    $b < \phi a$) and they are all higher than Bob's, then Alice has a
    winning strategy.
  \end{itemize}
\end{mainlemma}

We then use this lemma to prove our main result:

\begin{definition*}
  Say that a result holds \emph{generically} if it holds with
  probability approaching one as the number of cards, $n$, goes to
  infinity.
\end{definition*}

\begin{maintheorem}
  In the original game with unbiased dealing, generically neither
  player has a winning strategy.  Moreover, assume the cards are dealt
  randomly so that Alice receives $r \ge 1$ cards for every card of Bob.
  If $r < \phi$, generically neither player has a winning strategy,
  while if $r > \phi$, generically Alice has a winning strategy.
\end{maintheorem}

This result determines the limiting probability that one of the players
has a winning strategy, an open problem posed by Peter Winkler in his
book ``Mathematical Mind-Benders'', where it is attributed to the first
author.~\cite{Winkler}

\section*{Proofs}

One interesting quirk in Peer Pressure is that it is not immediately
obvious that having better cards is necessarily advantageous.  Of
course, a better card will be more likely to win in any given round, but
the card may also end up in the hands of the opponent, who may then use
it to his advantage.  We begin by proving that this is not a problem.
This lemma is not essential to the later results, but it does simplify
their proofs.

\begin{definition*}
  We say a collection of cards $C$ is \emph{at least as good as} another
  collection $C'$ if for all positive integers $k$, either the $k$th
  highest card in $C$ is at least as high as the $k$th highest card in
  $C'$ or $C'$ has less than $k$ cards.
\end{definition*}
\begin{monotonicitylemma}
  Having better cards doesn't hurt.  That is, if in a certain position
  Alice's cards are replaced with better cards and/or Bob's cards are
  replaced with worse cards, then Alice is no less likely to win.  In
  particular, if Alice had a winning strategy before the replacement,
  she still does afterward.
\end{monotonicitylemma}
\begin{proof} We prove the result step by step.

  \strong{Extra cards don't hurt.}  Suppose Alice receives extra cards,
  but no other change occurs.  Then clearly she is no worse off because
  she can play the same strategy as before, ignoring her extra cards.
  If she won before, she still wins.  (If she lost before, she now has
  extra cards that may or may not help in the end.)

  \strong{Losing cards doesn't help.}  By symmetry, if Bob has cards
  taken away, but no other change occurs, he is no better off.

  \strong{Receiving a card from the opponent doesn't hurt.}  Suppose Bob
  gives a card to Alice, but no other change occurs.  Then Alice is no
  worse off because this is equivalent to Alice gaining a card and Bob
  losing a card.

  \strong{Slightly improving one card doesn't hurt.}  Suppose Alice has
  a card $A$ and Bob has a card $B$ such that $A<B$ before the
  replacement and $A>B$ after the replacement, but no other change
  occurs.  (In particular, there are no cards of rank between $A$ and
  $B$.)  Then let Alice play exactly as before \emph{until} one of these
  cards is played by either player.  If both cards are played
  simultaneously, then Alice is no worse off because she wins a card
  instead of Bob.  Indeed, this is equivalent to Bob winning the battle
  (as before the replacement), followed by Bob giving Alice a card.  If
  only one card is played and it wins, then it is removed from play and
  so the relative ranking of $A$ and $B$ doesn't matter anyway.  If only
  one card is played and it loses, then one of the players will end up
  holding both $A$ and $B$.  Again, the relative ranking doesn't matter
  because Alice can pretend to switch the cards.

  \strong{Having better cards doesn't hurt.}  The general case consists
  of performing the above modifications one by one.  This may be
  accomplished, for example, by first improving Alice's best card, then
  her next best, and so on, and afterward, giving Alice extra cards and
  taking cards away from Bob.
\end{proof}

We may now prove our main lemma.

\begin{mainlemma}[redux]
  Let $\phi = \frac{1+\sqrt{5}}{2} \approx 1.61803$ be the golden ratio,
  which notably satisfies $1 + \phi = \phi^2$.
  \begin{itemize}
  \item
    If Alice has more than $\phi$ times as many cards as Bob (that is,
    $a > \phi b$), then Alice has a winning strategy.
  \item
    If Alice has more than $1/\phi$ times as many cards as Bob (that is,
    $b < \phi a$) and they are all higher than Bob's, then Alice has a
    winning strategy.
  \end{itemize}
\end{mainlemma}
\begin{proof}
  We prove the result by induction on the total number of cards, $a +
  b$.  Note that both results are certainly true when $a = 0$ or $b =
  0$.

  \strong{Alice has many cards.}  Suppose that $a > \phi b > 0$.  By the
  Monotonicity Lemma, we may assume that Bob has the $b$ highest cards,
  since this is the worst possible situation for Alice.  In order to
  win, Alice plays her lowest current card until all of Bob's cards are
  less than hers.  In any battle, Bob may either lose one of his high
  cards and receive one of Alice's low cards or Bob may give Alice back
  one of her original low cards.  In other words, Alice loses a card if
  and only if one of Bob's high cards is discarded.  Therefore,
  afterward, Bob has at most $b$ cards, all of which are lower than
  Alice's $a-b$ cards.  The result holds by induction because $b <
  \phi(a-b)$.

  \strong{Alice has enough high cards.}  Suppose that $\phi a > b > 0$
  and all of Alice's cards are higher than Bob's.  In order to win,
  Alice plays each of her $a$ cards once.  She will win every battle, so
  afterward she will still have $a$ cards, while Bob will have $b-a$
  cards.  The result holds by induction because $a > \phi (b-a)$.
\end{proof}

Armed with the Main Lemma, we prove the Main Theorem.  Recall that a
result holds \emph{generically} if it holds with probability approaching
one as the number of cards, $n$, goes to infinity.

\begin{maintheorem}[redux]
  In the original game with unbiased dealing, generically neither
  player has a winning strategy.  Moreover, assume the cards are dealt
  randomly so that Alice receives $r \ge 1$ cards for every card of Bob.
  If $r < \phi$, generically neither player has a winning strategy,
  while if $r > \phi$, generically Alice has a winning strategy.
\end{maintheorem}
\begin{proof}
  First of all, because only the relative ordering of the cards matters,
  we assume the cards are numbered $1$ to $n$.  Also, recall that if
  Alice has a mixed (probabilistic) winning strategy, then she also has
  a pure (deterministic) winning strategy.

  We begin with the weaker result with unbiased dealing.  Divide the
  cards into five equally-sized intervals: $C_1=\halfopen{0}{n/5}$,
  $C_2=\halfopen{n/5}{2n/5}$, $\dotsc$, $C_5 = \halfopen{4n/5}{n}$.
  Suppose Alice reveals her pure strategy in advance to Bob.  We will
  show how Bob can use these intervals to defeat it.  Specifically, for
  $1 \le i < 5$, he will use his cards in $C_{i+1}$ to defeat Alice's
  cards in $C_i$.  Finally, he will use his leftover cards to defeat
  Alice's~$C_5$.  See the Figure for a visual explanation.

  \makeatletter \renewcommand{\fnum@figure}{\figurename} \makeatother
  \begin{figure}[ht]
    \begin{center}
      \begin{pspicture}(-2,-0.4)(10,2)

        \psline[linestyle=dotted]( 0,0)( 0,2)
        \psline[linestyle=dotted]( 2,0)( 2,2)
        \psline[linestyle=dotted]( 4,0)( 4,2)
        \psline[linestyle=dotted]( 6,0)( 6,2)
        \psline[linestyle=dotted]( 8,0)( 8,2)
        \psline[linestyle=dotted](10,0)(10,2)

        \uput[270](1,0){$C_1$}
        \uput[270](3,0){$C_2$}
        \uput[270](5,0){$C_3$}
        \uput[270](7,0){$C_4$}
        \uput[270](9,0){$C_5$}

        \rput(-1,0.25){\strong{Alice}}
        \pspolygon( 0,0)( 0,0.5)( 2,0.5)( 2,0)
        \pspolygon( 2,0)( 2,0.5)( 4,0.5)( 4,0)
        \pspolygon( 4,0)( 4,0.5)( 6,0.5)( 6,0)
        \pspolygon( 6,0)( 6,0.5)( 8,0.5)( 8,0)
        \pspolygon( 8,0)( 8,0.5)(10,0.5)(10,0)

        \rput(-1,1.75){\strong{Bob}}
        \pspolygon( 0,1.5)( 0,2)( 2,2)( 2,1.5)
        \pspolygon( 2,1.5)( 2,2)( 4,2)( 4,1.5)
        \pspolygon( 4,1.5)( 4,2)( 6,2)( 6,1.5)
        \pspolygon( 6,1.5)( 6,2)( 8,2)( 8,1.5)
        \pspolygon( 8,1.5)( 8,2)(10,2)(10,1.5)

        \pspolygon[fillstyle=vlines]( 0,1.5)( 0,2)( 2,2)( 2,1.5)
        \pspolygon[fillstyle=vlines](3.26,1.5)(3.26,2)( 4,2)( 4,1.5)
        \pspolygon[fillstyle=vlines](5.26,1.5)(5.26,2)( 6,2)( 6,1.5)
        \pspolygon[fillstyle=vlines](7.26,1.5)(7.26,2)( 8,2)( 8,1.5)
        \pspolygon[fillstyle=vlines](9.26,1.5)(9.26,2)(10,2)(10,1.5)
        \pspolygon[fillstyle=vlines]( 8,0)( 8,0.5)(10,0.5)(10,0)

        \psline[arrows=->](2.63,1.4)(1,0.6)
        \psline[arrows=->](4.63,1.4)(3,0.6)
        \psline[arrows=->](6.63,1.4)(5,0.6)
        \psline[arrows=->](8.63,1.4)(7,0.6)
      \end{pspicture}
    \end{center}
    \caption{The unshaded intervals illustrate how Bob uses slightly
      more than a $1/\phi$-proportion of his cards in $C_{i+1}$ to
      defeat Alice's cards in $C_i$.  Bob's leftover cards in all of his
      $C_i$, represented by shaded intervals, are sufficient in number
      to overwhelm Alice's cards in $C_5$.}
  \end{figure}
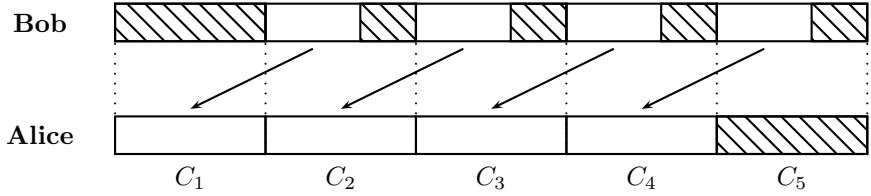

  We expect each player to receive half of the cards in each interval,
  so by the law of large numbers, Bob will generically receive at least
  $.099n$ cards from each interval; Alice will receive at most the
  remainder, $.101n$.  Now suppose that only Bob's cards in $C_2$ and
  Alice's cards in $C_1$ are under consideration.  By the Main Lemma,
  Bob may use $.063n > .101n/\phi$ cards from $C_2$ to defeat Alice's
  cards in $C_1$ (note that all of the cards in $C_2$ are higher than
  those in $C_1$), leaving at least $(0.099 - 0.063)n = 0.036n$ unused
  cards left over.  If he does similarly for his $C_3$ through $C_5$,
  Bob will have at least $4(.036n) + .099n = .243n$ cards left over.
  Again by the Main Lemma, Bob may use these cards to defeat Alice's
  $C_5$ because $.243 / .101 > \phi$.

  In the previous paragraph, we pretended that Bob may consider the game
  as the sum of five independent games.  However, this is justified
  because Alice has revealed her pure strategy in advance.  Because Bob
  knows where Alice will play, he may use the appropriate cards to
  defeat her.  Note that Bob may choose \emph{beforehand} which cards
  are allocated where, so it does not matter in what order Alice plays;
  in particular, Bob can choose his ``leftover cards'' beforehand, as
  all that matters is their number.  Therefore generically, Alice has no
  winning strategy and by symmetry, neither does Bob.

  Now we prove the stronger result with dealing biased towards Alice.
  If $r > \phi$, this is easy.  Generically, Alice will have more than
  $\phi$ times as many cards as Bob and thus win by the Main Lemma.

  Now suppose that $r < \phi$ is fixed.  We follow the same approach as
  before.  Divide the cards into $k$ equally-sized intervals $C_i =
  \halfopen{\frac{i-1}{k}n}{\frac{i}{k}n}$, where $k$ will be chosen
  later to depend only on $r$ and not on $n$.  In each interval, we
  expect players to receive cards in an $r : 1$ proportion.  By the law
  of large numbers, for any constant $\delta > 0$, Bob generically
  receives at least $1-\delta$ of the number of cards he expects in
  \emph{each} of the intervals.  (Note that we crucially use here that
  $k$ does not depend on $n$.)  By the Main Lemma, we may choose
  $\delta$ so small (but independent of $k$) that Bob may use his cards
  in $C_{i+1}$ to defeat Alice's cards in $C_i$ for all $1 \le i < k$
  and still have a positive proportion of his cards left over in each
  interval.

  Now choose $k$ so large (but independent of $n$) that Bob's remaining
  cards in $C_1$ and his leftover cards from all of the other $C_i$ are
  more than $\phi$ times the number of Alice's cards in $C_k$.  This is
  possible because we insured that each $C_i$ has at least a fixed
  positive proportion of cards left over, so Bob may overwhelm Alice
  with his extra cards.  Now, as before, Bob's leftover cards defeat
  Alice's $C_k$ by the Main Lemma.  (Again, Bob's cards may be allocated
  before any actual play.)  Finally, if Alice reveals her strategy in
  advance, Bob may combine his strategies on all of the intervals $C_i$
  to defeat her.

  Therefore, if $r < \phi$, Alice generically does not have a winning
  strategy.  Bob generically doesn't have a winning strategy either,
  because his cards are even worse than in the unbiased case.
\end{proof}

\section*{Further directions}

The results in this paper may be continued in a few natural directions.

For example, by using techniques similar to those presented above, Jacob
Fox has determined the threshold for the number of battles a player can
guarantee winning in the unbiased model.~\cite{Fox} In particular, if
$f(n)$ is a function that grows slower than $\sqrt{n}$ (using Landau's
asymptotic notation, $o(\sqrt{n})$), then generically both of the
players may guarantee winning at least $f(n)$ of the battles.  However,
if $f(n)$ is a function that grows faster than $\sqrt{n}$ (using
Landau's asymptotic notation, $\omega(\sqrt{n})$), then generically
neither player may guarantee winning at least $f(n)$ of the battles.
The threshold $\sqrt{n}$ comes from the central limit theorem.

In another direction, note that the Main Lemma classifies some of the
hands where Alice has a winning strategy, and it can also be used to
classify hands where neither player has a winning strategy (as in the
Main Theorem).  We leave as an open problem whether or not one may prove
stronger results about winning strategies:

\begin{problem*}
  Classify the situations when a given player has a winning strategy.
\end{problem*}

\section*{Acknowledgments}

The authors wish to thank Dan Cranston and Jacob Fox for helpful
discussions.

\bibliographystyle{amsalpha}

\begin{thebibliography}{Win07}

\bibitem[Fox08]{Fox}
Jacob Fox, personal communication, November 2008.

\bibitem[Win07]{Winkler}
Peter Winkler, \emph{Mathematical mind-benders}, ch.~11, p.~134, A K Peters
  Ltd., Wellesley, MA, 2007.

\end{thebibliography}
\providecommand{\bysame}{\leavevmode\hbox to3em{\hrulefill}\thinspace}
\providecommand{\MR}{\relax\ifhmode\unskip\space\fi MR }
\providecommand{\MRhref}[2]{%
  \href{http://www.ams.org/mathscinet-getitem?mr=#1}{#2}
}
\providecommand{\href}[2]{#2}

\end{document}